\documentclass[12pt,reqno]{amsart}
\newtheorem{thm}{Theorem}[section]
\newtheorem{defi}{Definition}[section]

\newtheorem{cor}{Corollary}[section]
\newtheorem{pr}{Proposition}[section]
\theoremstyle{definition}
\newtheorem{rem}{Remark}[section]

\newcommand{\be}{\begin{equation}}
\newcommand{\ee}{\end{equation}}
\newcommand{\bea}{\begin{eqnarray}}
\newcommand{\eea}{\end{eqnarray}}
\newcommand{\beb}{\begin{eqnarray*}}
\newcommand{\eeb}{\end{eqnarray*}}
\usepackage{amssymb,amsfonts,amsthm,setspace,indentfirst,mathrsfs}
\usepackage{minibox}
\usepackage{enumitem}
\usepackage{cite}

\usepackage{matlab-prettifier}

\usepackage{hyperref}

\numberwithin{equation}{section}
\begin{document}
%
\title{ Ricci solitons and curvature inheritance on Robinson-Trautman spacetimes}
\author[A. A. Shaikh and B. R. Datta]{Absos Ali Shaikh$^{*1}$ and Biswa Ranjan Datta$^2$}
\date{\today}
\address{\noindent\newline$^{1,2}$ Department of Mathematics,
	\newline University of Burdwan, 
	\newline Golapbag, Burdwan-713104,
	\newline West Bengal, India} 
\email{aask2003@yahoo.co.in$^1$, aashaikh@math.buruniv.ac.in$^1$}
\email{biswaranjandatta2019@gmail.com$^2$}

%
%
\dedicatory{}
\begin{abstract}
The purpose of the article is to investigate the existence of Ricci solitons   and the nature of curvature inheritance as well as collineations on the Robinson-Trautman (briefly, RT) spacetime. It is shown that under certain conditions RT spacetime admits almost Ricci soliton, almost $\eta$-Ricci soliton, almost gradient $\eta$-Ricci soliton. As a generalization of curvature inheritance \cite{Duggal1992} and curvature collineation \cite{KLD1969}, in this paper, we introduce the notion of \textit{generalized curvature inheritance} and examine if RT spacetime admits such a notion. It is shown that RT spacetime also realizes the generalized curvature (resp. Ricci, Weyl conformal, concircular, conharmonic, Weyl projective) inheritance. Finally, several conditions are obtained, under which RT spacetime possesses curvature (resp. Ricci, conharmonic, Weyl projective) inheritance as well as curvature  (resp. Ricci, Weyl conformal, concircular, conharmonic, Weyl projective) collineation, and we have also introduced the concept of generalized Lie inheritance and showed that RT spacetime realizes such a notion. 

\end{abstract}
%
\subjclass[2020]{53C25, 53C80,53C50, 83C05, 83C20}
\keywords{Robinson-Trautman spacetime; Ricci soliton; almost Ricci soliton; almost $\eta$-Ricci soliton; curvature inheritance; Ricci inheritance; Lie inheritance; Ricci collineation; curvature collineation; Weyl conformal collineation; concircular collineation; conharmonic collineation; Weyl projective collineation}
\maketitle
%

\section{\bf Introduction}\label{intro}

During the investigation of compact $3$-dimensional manifolds with positive Ricci curvature, in 1982, Hamilton \cite{Hamilton1982} established the concept of Ricci flow, a process of evolving a Riemanninan metric over time. The main idea of Hamilton was to smooth out the singularities of the metric by establishing a new type of non-linear diffusion equation. Ricci solitons \cite{Hamilton1988} are the self-similar solutions of the Ricci flow, which is a natural generalization of Einstein metrics \cite{Bess87,Brink1925,S09}. There are several generalizations of Ricci solitons such as almost Ricci solitons, $\eta$-Ricci solitons, gradient $\eta$-Ricci solitons etc. During last three decades, lots of results (see, \cite{AliAhsan2013, AliAhsan2015, Ahsan2018}) on Ricci solitons have been appeared and now it is an active area of research in differential geometry.      \\

If the Ricci curvature $S$ of a Riemannian manifold $(\mathcal Q,g)$ satisfies 
 $$\frac{1}{2}\mathscr{L}_\xi g_{ij}+S_{ij}-\mu g_{ij}=0$$
with the constant $\mu$ and the Lie derivative $\mathscr{L}_\xi$ in the direction of the soliton vector field $\xi$, then $\mathcal Q$ is said to be a Ricci soliton. If $\xi=\nabla \zeta$ for a smooth function $\zeta$ on $\mathcal Q$, then $\mathcal Q$ is said to be gradient Ricci soliton with $\zeta$ as the potential function and it satisfies
$$S_{ij}+(\nabla^2f)_{ij}-\mu g_{ij}=0.$$ 
 The Ricci soliton is expanding, steady or shrinking  as per the conditions $\mu<0$, $\mu=0$ or $\mu>0$ respectively. If $\mu$ is non-constant, specially a smooth function, then it is called an almost Ricci soliton \cite{Pigola2011}. For a Killing soliton vector field $\xi$, the Ricci soliton becomes an Einstein manifold. If $(\mathcal Q,g)$ admits a non-zero $1$-form $\eta$ satisfying
$$\frac{1}{2}\mathscr{L}_\xi g_{ij}+S_{ij}-\mu g_{ij}+\lambda (\eta\otimes\eta)_{ij}=0,$$
$\mu,\lambda$ being constants, then it is called a $\eta$-Ricci soliton \cite{Cho2009}. Again, the soliton is known as an almost $\eta$-Ricci soliton \cite{Blaga2016} if $\mu, \lambda$ are allowed to be smooth functions. \\

Again, the geometrical symmetries perform a pivotal role in the theory of general relativity as realizing the Einstein's field equations (EFE), the establishment of gravitational potentials can be obtained by imposing the symmetries. Such geometrical symmetry leads to conservation laws in the appearance of the first integral of a dynamical system, which are recognized as Killing vector fields \cite{Noether_1918, Davis_Katzin_1962}. Several kinds of geometrical symmetries in spacetimes can be achieved if the Lie derivative of certain tensor with respect to some vector field vanishes. The vanishing Lie derivative of a geometric quantity with respect to some vector field usually represents the preservation of the geometric quantity in the direction of the vector field. Such symmetries are known as collineations and they can be described if in the direction of some vector field, certain geometric quantities such as metric tensor, Ricci tensor, stress-energy momentum tensor, Riemann curvature tensor, Weyl conformal curvature tensor, Weyl projective curvature tensor etc. remain invariant. The underlying vector field may be null, spacelike or timelike. 
The role of collineations in general relativity was demonstrated by Katzin et al. \cite{KLD1969,KLD1970}. 
In the context of gravitational radiation, investigation on conservation laws consists of a fundamental symmetry property of a spacetime, known as curvature collineation as the relations between groups of motions, collineations and conservation laws of energy, momentum are crucial in understanding the geometry of a spacetime. The notion of curvature collineation is investigated in different frames by Hall and Shabbir \cite{Hall_2001}, Singh and Sharma \cite{Singh_1975}, Tiwari \cite{Tiwari_2005} with their relevance in general relativity. Also, Ahsan \cite{Ahsan1978}, Ahsan and Husain \cite{AH1980} studied several kinds of collineations to find the relations between physical properties and geometrical symmetries of elecrtomagnatic fields and showed that motion and Maxwell collineation are not comparable for null electromagnatic field. Recently, Haesen and Verstraelen \cite{Haesen_2007} introduced pseudosymmetric collineation and studied it in vacuum pp-wave spacetime along with some interesting  relations with other types of transformations. Using the Nijenhuis tensor, the notion of torsion collineation was introduced by Ahsan \cite{Ahsan1987}. Also, several geometrical aspects of collineation have been further studied by Ahsan et al. \cite{Ahasan2005, Ahsan1995, Ahsan1996, Ahsan1977_1055,Ahsan1977_231,AhsanAli2014}, Ali and Ahsan \cite{AliAhsan2012,AA2012}. Again, Duggal \cite{Duggal1992} established the generalized concept of symmetry, called curvature inheritance, extending the notion of curvature collineation. The notion of curvature inheritance have been studied by Coley and Tupper \cite{Coley_1989} and Shaikh et al. \cite{Shaikh_Ali_Salman_2023}. The detailed reviews of various kinds of symmetry inheritance can be found in \cite{Duggal_1999, Hall_1991, Duggal_1987}. During last two decades such notion of symmetries have become an increasingly prevalent area of research in general relativity.\\


The Riemannian manifold $(\mathcal Q,g)$ admits curvature (resp., Ricci, Weyl conformal, concircular, conharmonic, Weyl projective) collineation if the Lie derivative $\mathscr L_\xi R$ of the Riemann curvature tensor $R$ (resp., Ricci curvature tensor $S$, Weyl conformal curvature tensor $C$, concircular curvature tensor $W$,  conharmonic curvature tensor $K$, Weyl projective curvature tensor $P$) vanishes, where $\xi$ is a smooth vector field. Generalizing the notion of curvature collineation, Duggal \cite{Duggal1992} introduced the concept of curvature inheritance. In the present paper, we have introduced the notion of \textit{generalized curvature inheritance} (see Definition \ref{GI}), which includes the notion of curvature inheritance as well as curvature collineation. Then, we have showed that RT spacetime fulfills generalized curvature (resp. Ricci, Weyl conformal, concircular, conharmonic, Weyl projective) inheritance. Also, we have deduced several conditions under which RT spacetime realizes the curvature (resp. Ricci, conharmonic, Weyl projective) inheritance and curvature (resp. Ricci, Weyl conformal, concircular, conharmonic, Weyl projective) collineation.\\

RT spacetime \cite{RT1960} was introduced in 1960 as a class of algebraically special solutions of EFE having a repeated principal null direction associated with a shearless, expanding and non-twisting null geodesic congruence. During the collision of two black holes the estimation of mass loss in the final phase can be derived by the RT spacetime. Also, the RT spacetime can represent the model of gravitational radiation outgoing from spatially bounded sources. Describing topological equivalent two spheres instead of strictly two spheres, the RT spacetimes generalize the notion of spherical symmetry. The RT spacetime is a crucial class of exact solutions of EFE as it is used to describe several models of gravitational waves, black holes and cosmology. The Petrov type II RT spacetime is given as follows (\cite{RT1960,RT1962,Podolsky,SKMHH2003}):
$$ds^2=-2(U^0-2\gamma^0r-\Psi_2^0r^{-1})du^2+2dudr-\frac{r^2}{2\Omega^2}d\zeta d\bar\zeta$$
where $U^0, \gamma^0, \Psi_2^0$ are constants and $\Omega$ is a everywhere non-vanishing function of $\zeta, \bar\zeta$. Setting $\zeta=x+iy$, the RT metric can be expressed as follows:
\begin{equation}\label{RTM}
ds^2=-2\left(a-2br-\frac{d}{r}\right)dt^2+2dtdr-\frac{r^2}{f^2}(dx^2+dy^2),
\end{equation}
where $a,b,d$ are constants and $f$ is a nowhere vanishing function of $x$ and $y$. Again, for the warping function $r$, the  metric tensor $g$ of RT spacetime is the warped product $\bar{g}\times_r \tilde{g}$ with  base metric $\bar g$ and fiber metric $\tilde g$ given by 

\begin{equation}\nonumber
	\bar g=\begin{pmatrix}
		-2(a-2br-\frac{d}{r}) & 1\\
		1 & 0
	\end{pmatrix}\ \  \text{ and } \ \ 
\tilde g = \begin{pmatrix}
	-\frac{1}{f^2} & 0\\
	0 & -\frac{1}{f^2}
\end{pmatrix}.
\end{equation}
\noindent In 2018, Shaikh et al. \cite{SAA18} explored several curvature properties of RT spacetime and showed that RT spacetime is  2-quasi-Einstein, Roter type manifold and it admits Weyl pseudosymmetric, Ricci pseudosymmetric, pseudosymmetric Weyl conformal curvature tensor, Riemann compatible Ricci tensor, recurrent Weyl conformal curvature 2-forms, pseudosymmetric energy momentum tensor. Numerous  relativistic theories on different spacetimes can be found in \cite{warsaw_1962}. Recently, the curvature properties of several spacetimes such as Melvin magnatic spacetime \cite{SAAC20}, Nariai spacetime \cite{SAAC20N}, Lemaitre-Tolman-Bondi spacetime \cite{SAACD22}, generalized pp-wave spacetime \cite{SBH21}, Kantowski-Sachs spacetime \cite{SC21}, Vaidya-Bonner spacetime \cite{SDC21}, Siklos spacetime \cite{SDKC19}, Som-Raychaudhuri spacetime \cite{SK16srs}, Lifshitz spacetime \cite{SSC19}, interior black hole spacetime \cite{SDHK20}, Vaidya spacetime \cite{SKS2019}, $(t-z)$-type plane wave spacetime \cite{Eyasmin2021t-z}, Morris-Thorne wormhole spacetime \cite{Eyasmin2022Morris} etc. have been investigated.  \\


Following the resolution of Hamilton's theorized  conjuncture, the issue of finding the generalization of Ricci soliton has gained significant importance. The notions of almost Ricci soliton, $\eta$-Ricci soliton, almost $\eta$-Ricci soliton are quite noteworthy among such generalizations. General relativity encountered several applications of almost Ricci solitons, as demonstrated by Duggal \cite{Duggal_2017}. The different aspects of Ricci solitons and its generalizations have been recently studied by  Ahsan et al. \cite{Ahsan2018,AliAhsan2013, AliAhsan2015},  Blaga \cite{Blaga2016} and Deshmukh et al. \cite{Deshmukh_2011, Deshmukh_2019, Deshmukh_2020}, Shaikh et al. \cite{SMM2022,SCM2021, SCM2021yamabe, SDAA2021, SMM2021, SM2021}. But the nature of Ricci solitons on RT spacetime remains to be investigated yet. Therefore, a natural question arises: whether RT spacetime admits Ricci solitons or not? If so, then what would be the nature of such solitons or under which conditions it admits Ricci solitons and  collineations? The answer of these questions have been given in the present paper. The main objective of this article is to investigate whether the RT spacetime admits curvature inheritance, Ricci solitons and collineations with different curvature tensors such as Ricci, conharmonic, projective curvature tensor under certain conditions.  The technique of the present article is unconventional as we have calculated various components of different related tensors of Ricci solitons and collineations from the RT metric, and then we have checked \& verified the concerned governing equations with the help of a program developed by the first author in Wolfram Mathematica, which is explored elegantly in entire Section 6.\\

The outline of the paper is delineated in the following ways: Section 2 is devoted to the basic rudiments of the geometric and curvature properties of RT spacetimes. Section 3 deals with the investigation of several Ricci solitons on RT spacetime and their nature. It is shown that under various conditions RT spacetime can be exhibited as an example of almost Ricci soliton, almost $\eta$-Ricci soliton, almost gradient $\eta$-Ricci soliton.  Section 4 is engaged with the investigation of the inheritance and collineations with respect to Riemann (resp. Ricci, Weyl conformal, projective, concircular and conharmonic) curvature tensor. The conclusion of this article along with a discussion of potential future research have been included in Section 5. Finally, the last section ``Appendix" is solely devoted to the methodology, which provides a gradual analysis of the Wolfarm Mathematica program with the implementation techniques of the program codes in the present study.

\section{Preliminaries}
Let $(\mathcal Q,g)$ be a smooth connected $n$-dimensional $(n\ge 3)$ semi-Riemannian manifold with the Lie algebra $\chi(\mathcal Q)$ of all smooth vector fields on $\mathcal Q$. Now, the Kulkarni-Nomizu product $\phi \wedge \psi$ of two $(0,2)$ type symmetric tensors $\phi$ and $\psi$ is defined as (\cite{Glog02}): 
\bea
(\phi \wedge \psi)_{pqrs}= \phi_{ps}\psi_{qr} - \phi_{pr}\psi_{qs}+ \phi_{qr}\psi_{ps} - \phi_{qs}\psi_{pr}.\nonumber
\eea
Let $\mathcal T^r_s(\mathcal Q)$ ($r,s\geq 0$ are non-negative integers) be the space of tensor fields of type $(r,s)$ on the manifold $\mathcal Q$. A tensor $\Pi\in \mathcal T^1_3(\mathcal Q)$ is a generalized curvature tensor  if $\Pi$ satisfies the following (\cite{Nomizu1971,SDHJK15}):
	\begin{enumerate}[label=(\roman*)]
		\item $\Pi(\mathcal O_1, \mathcal O_2)\mathcal O_3 + \Pi(\mathcal O_2, \mathcal O_1)\mathcal O_3=0$,
		
		\item $\Pi(\mathcal O_1, \mathcal O_2,\mathcal O_3,\mathcal O_4)=\Pi(\mathcal O_3, \mathcal O_4,\mathcal O_1,\mathcal O_2)$ and
		
		\item $\Pi(\mathcal O_1, \mathcal O_2)\mathcal O_3 + \Pi(\mathcal O_2, \mathcal O_3)\mathcal O_1 + \Pi(\mathcal O_3, \mathcal O_1)\mathcal O_2 = 0$ (Bianchi's $1^{st}$ identity),
	\end{enumerate}
where, $\Pi(\mathcal O_1, \mathcal O_2,\mathcal O_3,\mathcal O_4)=g(\Pi(\mathcal O_1, \mathcal O_2)\mathcal O_3,\mathcal O_4)$ and for $i=1,2,3,4$ $\mathcal O_i\in\chi(\mathcal Q)$. Further, if $\Pi$ satisfies Bianchi's  $2^{nd}$ identity, i.e.,
 
$$(\nabla_{\mathcal O_1}\Pi)(\mathcal O_2, \mathcal O_3)\mathcal O_4 + (\nabla_{\mathcal O_2}\Pi)(\mathcal O_3, \mathcal O_1)\mathcal O_4 + (\nabla_{\mathcal O_3}\Pi)(\mathcal O_1, \mathcal O_2)\mathcal O_4 = 0,$$
$\nabla$ being the Levi-Civita connection on $\mathcal Q$, then $\Pi$ is said to be a proper generalized curvature tensor. Some crucial examples of generalized curvature tensors are the Riemann curvature $R$, Gaussian curvature $G$, conharmonic curvature $K$, concircular curvature $W$, Weyl conformal curvature $C$, which are defined respectively as follows (\cite{SAR13,SB08,SK14,SKppsn}):

\bea
R_{pqrs} &=& g_{p\alpha}(\partial_{s}\Gamma^{\alpha}_{qr} - \partial_r \Gamma^{\alpha}_{qs} + \Gamma^{\beta}_{qr}\Gamma^{\alpha}_{\beta s} - \Gamma^{\beta}_{qs}\Gamma^{\alpha}_{\beta r}), \nonumber \\
G_{pqrs} &=& \frac{1}{2}\left( g \wedge g\right)_{pqrs}, \nonumber \\
K_{pqrs} &=&  R_{pqrs} - n_3(S\wedge g)_{pqrs}, \nonumber\\
W_{pqrs} &=&  R_{pqrs} - \kappa n_1 n_2 G_{pqrs}  \nonumber \\
C_{pqrs} &=&  R_{pqrs} - n_3 (S\wedge g)_{pqrs} + \kappa n_2 n_3 G_{pqrs} \nonumber
\eea
where $n_\lambda=\frac{1}{1+n-\lambda}$ ($\lambda=1,2,3$), $\kappa$ is the scalar curvature, $\Gamma^p_{qr}$ are connection coefficients and $\partial_u=\frac{\partial}{\partial x^u}$.
But, the Weyl projective curvature defined by
$$P^u_{pqr}= R^u_{pqr} -2 n_2\delta^u_{[p}S_{q]r},$$
is not a generalized curvature tensor. We note that $W$, $C$, $K$ are not proper generalized curvature tensor.

\begin{defi}
	If $\mathscr T\in \mathcal T^0_4(\mathcal Q)$, then the manifold $\mathcal Q$ is said to possess $\mathscr T$ collineation for some $\xi\in \chi(\mathcal Q)$ if 
$$\mathscr{L}_\xi \mathscr T=0,$$
where $\mathscr L_\xi$ denotes the Lie derivative along the vector field $\xi$. If $\mathscr T=R$ (resp. $C, W, K, P$), then the  $\mathscr T$ collineation is said to be curvature (resp.  Weyl conformal,  concircular, conharmonic, Weyl projective) collineation.
\end{defi}

Generalizing the concept of curvature collineation, in 1992, Duggal \cite{Duggal1992} introduced the notion of curvature inheritance defined as follows:
$$	\mathscr L_{\xi}R=\lambda R,$$
where $\lambda$ is a smooth function and $\xi\in\chi(\mathcal Q)$. Now, generalizing the notion of curvature inheritance we introduce the concept of generalized $\mathscr T$ inheritance, which is defined as follows:
\begin{defi}\label{GI}
	If $\mathscr T\in \mathcal T^0_4(\mathcal Q)$, then the manifold $\mathcal Q$ is said to admit \textbf{generalized $\mathscr T$ inheritance} with respect to some vector field $\xi$ if 
	$$	\mathscr L_{\xi}\mathscr T=\lambda \mathscr T+\lambda_1 g\wedge g+\lambda_2 g\wedge S,$$
where $\lambda,\lambda_1,\lambda_2$ are smooth functions and $\wedge$ is the Kulkarni-Nomizu product. If $\mathscr T=R$ (resp. $C, W, K, P$), then the generalized $\mathscr T$ inheritance is said to be generalized curvature (resp.  Weyl conformal,  concircular, conharmonic, Weyl projective) inheritance.
\end{defi}
 If $\lambda_1,\lambda_2=0$ (resp.  $\lambda_1,\lambda_2,\lambda=0$), then the generalized curvature inheritance turns into curvature inheritance (resp. curvature collineation). 
 
 \begin{defi}
 	For  $\mathscr B\in \mathcal T^0_2(\mathcal Q)$, the manifold $\mathcal Q$ is said to be admitted $\mathscr B$-collineation if the Lie derivative of $\mathscr B$ with respect to some vector field $\xi$ vanishes, i.e.,
 	$$\mathscr{L}_\xi \mathscr B=0.$$
 \end{defi}
 If the tensor $\mathscr B$ is chosen as the metric tensor $g$, Ricci curvature tensor $S$ and stress-energy momentum tensor $T$, then the $\mathscr B$-collineations are also known as motion, Ricci collineation \cite{BQ1993} and matter collineation respectively. The underlying vector field for a motion is called Killing vector field. Generalizing the concept of Ricci collineation and Ricci (resp. Lie) inheritance, we introduced the notion of generalized Ricci (resp. Lie) inheritance, defined as follows:
\begin{defi}
The manifold $\mathcal Q$ is said to be admitted \textbf{generalized Ricci inheritance} if
 	$$\mathscr{L}_\xi S=\lambda_S S + \lambda_g g,$$
where $\lambda_S,\lambda_g$ are smooth functions. If $\lambda_g=0$, then it turns out to be Ricci inheritance and if $\lambda_S=0=\lambda_g$, then it is called Ricci collineation.
\end{defi} 

In particular, if the manifold $\mathcal Q$ is Einstein $\left(S=\frac{\kappa}{n}g \right)$ and realizes generalized Ricci inheritance, then the vector field $\xi$ is a conformal Killing vector field, i.e., the motion with respect to $\xi$ is a conformal motion. Moreover, in this case, if $\lambda_S$, $\lambda_g$ are constants, then $\xi$ is homothetic, and also, if the scalar curvature $\kappa=-n\frac{\lambda_g}{\lambda_S}$, then $\xi$ is a Killing vector field. 

\begin{defi}
The manifold $\mathcal Q$ is said to be admitted \textbf{generalized Lie inheritance} if
 	$$\mathscr{L}_\xi T=\lambda_T T + \lambda_g g,$$
where $\lambda_T,\lambda_g$ are smooth functions. If $\lambda_g=0$, then it turns out to be Lie inheritance and if $\lambda_T=0=\lambda_g$, then it is known as Lie symmetry along $\xi$.
\end{defi} 

We note that if the manifold is a spacetime and admits generalized Lie inheritance, then from EFE we conclude that the spacetime possesses the following relation:
$$\mathscr{L}_\xi S+a \mathscr{L}_\xi g+bg=\lambda_T S,$$
where $a=\Lambda-\frac{\kappa}{2}$ and $b=\Lambda-\frac{d\kappa(\xi)}{2}-K\lambda_g-\lambda_T\left(\Lambda-\frac{\kappa}{2} \right)$ with the cosmological constant $\Lambda$ and Einstein's gravitational constant $K$. If the manifold $\mathcal Q$ is Einstein, then the previous relation yields that $\xi$ is a conformal Killing vector field, i.e., the motion with respect to $\xi$ is the conformal motion.

The components of Riemann curvature tensor and  Ricci tensor of the RT spacetime are given as follows:

$$\begin{array}{c}
	R_{1212}=-\frac{2d}{r^3}, R_{1313}=\frac{2(d-2br^2)(d+r(-a+2br))}{r^2f^2}= R_{1424}, R_{1323}= \frac{d-2br^2}{rf^2}=R_{1424}, \\
	R_{3434}= \frac{r(2(d+r(-a+2br))+r(f_x^2+f_y^2-f(f_{xx}+f_{yy})))}{f^4};
\end{array}$$

 
$$\begin{array}{c}
	S_{11}= -\frac{8b(d+r(-a+2br))}{r^2}, S_{12}= -\frac{4b}{r}, 
	S_{33}= \frac{-2a+8br+f_x^2+f_y^2-f(f_{xx}+f_{yy})}{f^2}=S_{44}
\end{array}$$
The scalar curvature $\kappa$ is given by
$\kappa= -\frac{2}{r^2} \left[f_x^2+f_y^2-f(f_{xx}+f_{yy})-2a+12br\right]$.

\section{\bf Nature of Ricci solitons on Robinson-Trautman spacetime}\label{solitons}

For an $n$-dimensional smooth semi-Riemannian manifold $\mathcal Q$, the set $\mathcal S(\mathcal Q)$ of all Killing vector fields configures a Lie subalgebra of $\chi(\mathcal Q)$. The maximum number of linearly independent Killing vector fields in $\mathcal S(\mathcal Q)$ is less than or equal to $n(n+1)\slash 2$. The equality holds if $\mathcal Q$ is of constant scalar curvature. We note that the scalar curvature $\kappa$ of RT spacetime is not a constant. In RT spacetime, the vector field $\frac{\partial}{\partial t}$ is a Killing vector field, i.e.,  $\mathscr L_{\frac{\partial}{\partial t}} g$ vanishes. In this section, we have considered several non-Killing vector fields such as $\frac{\partial}{\partial r}$, $\frac{\partial}{\partial x}$, $\frac{\partial}{\partial y}$, $\mu_1\frac{\partial}{\partial r}+\mu_2\frac{\partial}{\partial x}+\mu_3\frac{\partial}{\partial y}$ and $\nabla r^2$ on the RT spacetime, where $\mu_1$, $\mu_2$, $\mu_3$ are constants.\\

For the non-Killing vector field $X_2=\frac{\partial}{\partial r}$, RT spacetime realizes the following relation:
\begin{eqnarray*}\label{e2}
\mathscr L_{X_2}g+\frac{2r}{F_1-2a+4br}S+\frac{8b}{F_1-2a+4br}g-\frac{2(-d+2br^2)}{r^2}\eta_1\otimes\eta_1=0,
\end{eqnarray*}
where $\eta_1=(1,0,0,0)$ is an $1$-form and $F_1=f_x^2+f_y^2-f(f_{xx}+f_{yy})$ is a smooth function such that $F_1-2a+4br\neq 0$. Thus, we can state the following:

\begin{pr}
	If $f_x^2+f_y^2-f(f_{xx}+f_{yy})-2a+4br\neq0$, then the RT spacetime admits an almost $\eta$-Ricci soliton given by 
	\begin{eqnarray*}
		\frac{1}{2}\mathscr L_{\xi}g+S+\frac{4b}{r}g-\frac{-d+2br^2}{r^2}\eta\otimes\eta=0,
	\end{eqnarray*}
where $\eta=(1,0,0,0)$ is an $1$-form and $\xi=\frac{\partial}{\partial r}$ is the soliton vector field. 
\end{pr}

Again, for the non-Killing vector field $X_3=\frac{\partial}{\partial x}$, RT spacetime fulfills the following relation:
\begin{eqnarray*}
	\mathscr L_{X_3}g-\frac{2r^2f_x}{f(2a-4br-F_1)}S-\frac{8brf_x}{f(2a-4br-F_1)}g=0,\label{e3}
\end{eqnarray*}
where $F_1=f_x^2+f_y^2-f(f_{xx}+f_{yy})$ with $F_1-2a+4br\neq0$. 
This implies the following:

\begin{pr}
	If $r^2f_x+f\{2a-4br-f_x^2-f_y^2+f(f_{xx}+f_{yy})\}=0$ and $f_x^2+f_y^2-f(f_{xx}+f_{yy})-2a+4br\neq0$, then the RT spacetime admits an almost Ricci soliton given by
	\begin{eqnarray*}
		\frac{1}{2}\mathscr L_{\xi}g+S+\frac{4b}{r}g=0,
	\end{eqnarray*}
where, $\xi=\frac{\partial}{\partial x}$ is the soliton vector field. Further, if $ r\neq 0 $ is a constant $c_1$ and $f\{f(f_{xx}+f_{yy})-f_x^2-f_y^2+2a-4bc_1\}+c_1^2f_x=0$, then  the RT spacetime possesses a Ricci soliton given by
	\begin{eqnarray*}
	\frac{1}{2}\mathscr L_{\xi}g+S+\frac{4b}{c_1}g=0.
\end{eqnarray*}
\end{pr}

Again, for the non-Killing vector field $X_4=\frac{\partial}{\partial y}$, RT spacetime satisfies the following relation:
\begin{eqnarray*}
	\mathscr L_{X_4}g-\frac{2r^2f_y}{f(F_1-2a+4br)}S-\frac{8brf_y}{f(F_1-2a+4br)}g=0,\label{e4}
\end{eqnarray*}
where $F_1=f_x^2+f_y^2-f(f_{xx}+f_{yy})$ and $F_1-2a+4br\neq0$. 
Thus we can state the following:

\begin{pr}
	If $r^2f_x+f\{-2a+4br+f_x^2+f_y^2-f(f_{xx}+f_{yy})\}=0$  and $f_x^2+f_y^2-f(f_{xx}+f_{yy})-2a+4br\neq0$, then the RT spacetime admits an almost Ricci soliton given by
	\begin{eqnarray*}
		\frac{1}{2}\mathscr L_{\xi}g+S+\frac{4b}{r}g=0,
	\end{eqnarray*}
where, $\xi=\frac{\partial}{\partial y}$ is the soliton vector field.
\end{pr}

Now, for the non-Killing vector field $V_1=\mu_1\frac{\partial}{\partial r}+\mu_2\frac{\partial}{\partial x}+\mu_3\frac{\partial}{\partial y}$ ($\mu_1,\mu_2,\mu_3$ are constants), the following relation  holds:
\begin{eqnarray*}
	\mathscr L_{V_1}g+\frac{2r\{-\mu_1f+r(\mu_2f_x+\mu_3f_y)\}}{f(2a-4br-F_1)}S+\frac{8b\{\mu_1f-r(\mu_2f_x+\mu_3f_y) \}}{f(2a-4br-F_1)}g-\frac{2\mu_1(-d+2br^2)}{r^2}\eta_2\otimes\eta_2=0,\label{e234}
\end{eqnarray*}
where $F_1=f_x^2+f_y^2-f(f_{xx}+f_{yy})$ with $F_1-2a+4br\neq0$ and $\eta_2=(1,0,0,0)$ is an $1$-form. This concludes the following:
\begin{pr}
	If $f\{2a-4br-f_x^2-f_y^2+f(f_{xx}+f_{yy})\}+r\{\mu_1f-r(\mu_2f_x+\mu_3f_y) \}=0$ and $f_x^2+f_y^2-f(f_{xx}+f_{yy})-2a+4br\neq0$, then the RT spacetime admits an almost $\eta$-Ricci soliton given by
	\begin{eqnarray*}
		\frac{1}{2}\mathscr L_{\xi}g+S-4bg-\frac{\mu_1(-d+2br^2)}{r^2}\eta\otimes\eta=0,
	\end{eqnarray*}
	where $\eta=(1,0,0,0)$ is an $1$-form and $\xi=\mu_1\frac{\partial}{\partial r}+\mu_2\frac{\partial}{\partial x}+\mu_3\frac{\partial}{\partial y}$ ($\mu_1,\mu_2,\mu_3$ are constants) is the soliton vector field. 
\end{pr}

Again, for the non-Killing vector field $X=grad(r^2)=2r\frac{\partial}{\partial t}-4\{d+r(-a+2br)\}\frac{\partial}{\partial r}$, the following relation holds:
\begin{eqnarray*}
	&&\mathscr L_Xg-\frac{4r\{3d+2r(-a+br)\}}{F_1-2a+4br}S-\frac{4\{2ad-16bdr+4abr^2+(-d+2br^2)F_1\}}{r(2a-4br-F_1)}g-4\eta_3\otimes\eta_3                                                           =0\label{gradr^2}
\end{eqnarray*}
where $\eta_3=(0,1,0,0)$ is an $1$-form and $F_1=f_x^2+f_y^2-f(f_{xx}-f_{yy})$ is a smooth function such that $F_1-2a+4br\neq0$. This leads to the following:

\begin{pr}
If $f_x^2+f_y^2-f(f_{xx}-f_{yy})+2\{3d+2r(-a+br)\}-2a+4br=0$ and $f_x^2+f_y^2-f(f_{xx}-f_{yy})-2a+4br\neq0$, then the RT spacetime admits an almost gradient  $\eta$-Ricci soliton given by  
\begin{eqnarray*}
	\frac{1}{2}\mathscr L_\xi g+S-\frac{2(-6bdr+4abr^2-5bdr^2+3d^2-2adr-4b^2r^3+4abr^3-2b^2r^4)}{r^2\{3d+2r(-a+br) \}}g-2\eta\otimes\eta=0
\end{eqnarray*}
where  $3d+2r(-a+br)\neq0$ with the $1$-form $\eta=(0,1,0,0)$ and $\xi=\nabla r^2$ is the soliton vector field.
\end{pr}

\begin{rem}
	If $\xi\in \chi(\mathcal Q)$ can be expressed as $\xi=\alpha\frac{\partial}{\partial r}+\beta\frac{\partial}{\partial x}+\gamma\frac{\partial}{\partial y}$ for some constants $\alpha, \beta, \gamma$, then the RT spacetime can not be an example of an $\eta$-Yamabe soliton or a Yamabe soliton for the soliton vector field $\xi$.
\end{rem}

\begin{thm}
The RT spacetime (\ref{RTM}) together with $f_x^2+f_y^2-f(f_{xx}+f_{yy})-2a+4br\neq0$ satisfies the following properties:
\begin{enumerate}[label=(\roman*)]
	\item for the soliton vector field $\xi=\frac{\partial}{\partial r}$, the RT spacetime admits an almost $\eta$-Ricci soliton  if $f_x^2+f_y^2-f(f_{xx}-f_{yy})-2a+4br-r=0$ with the $1$-form $\eta=(1,0,0,0)$, i.e.,
		\begin{eqnarray*}\nonumber
			\frac{1}{2}\mathscr L_{\xi}g+S+\frac{4b}{r}g-\frac{-d+2br^2}{r^2}\eta\otimes\eta=0,
	\end{eqnarray*}

	\item for the soliton vector field $\xi=\frac{\partial}{\partial x}$, the RT spacetime possesses an almost Ricci soliton  if $f\{f(f_{xx}+f_{yy})-f_x^2-f_y^2+2a-4bc_1\}+r^2f_x=0$, i.e.,
	\begin{eqnarray*}\nonumber
			\frac{1}{2}\mathscr L_{\xi}g+S+\frac{4b}{r}g=0,
		\end{eqnarray*}

	\item for the soliton vector field $\xi=\frac{\partial}{\partial y}$, the RT spacetime reveals an almost Ricci soliton  if $r^2f_x+f\{-2a+4br+f_x^2+f_y^2-f(f_{xx}+f_{yy})\}=0$, i.e.,
	\begin{eqnarray*}\nonumber
			\frac{1}{2}\mathscr L_{\xi}g+S+\frac{4b}{r}g=0,
		\end{eqnarray*}

	\item for the soliton vector field $\xi=\mu_1\frac{\partial}{\partial r}+\mu_2\frac{\partial}{\partial x}+\mu_3\frac{\partial}{\partial y}$, the RT spacetime admits an almost $\eta$-Ricci soliton  if $f\{2a-4br-f_x^2-f_y^2+f(f_{xx}+f_{yy})\}+r\{\mu_1f-r(\mu_2f_x+\mu_3f_y) \}=0$ with $\eta=(1,0,0,0)$, i.e.,
	\begin{eqnarray*}\nonumber
			\frac{1}{2}\mathscr L_{\xi}g+S-4bg-\frac{\mu_1(-d+2br^2)}{r^2}\eta\otimes\eta=0,
		\end{eqnarray*}

	\item for the soliton vector field $\xi=\nabla r^2$, the RT spacetime possesses an almost gradient  $\eta$-Ricci soliton 
	if $f_x^2+f_y^2-f(f_{xx}-f_{yy})+2\{3d+2r(-a+br)\}-2a+4br=0$ and $3d+2r(-a+br)\neq0$	with the $1$-form $\eta=(0,1,0,0)$, i.e., 
	\begin{eqnarray*}\nonumber
		\frac{1}{2}\mathscr L_\xi g+S-\frac{2(-6bdr+4abr^2-5bdr^2+3d^2-2adr-4b^2r^3+4abr^3-2b^2r^4)}{r^2\{3d+2r(-a+br) \}}g\\
		-2\eta\otimes\eta=0.
	\end{eqnarray*}

\end{enumerate}
\end{thm}

\section{\bf Curvature Inheritance and Collineations on Robinson-Trautman spacetime}\label{collineations}
If $f_x^2+f_y^2-f(f_{xx}+f_{yy})-2a+4br\neq 0$, then for the vector field $V=\mu_1 \frac{\partial}{\partial x}+\mu_2\frac{\partial}{\partial y}$, the following relation holds:
\begin{eqnarray}
	\mathscr L_{V}S=\lambda_S S+\lambda_g g,\notag
\end{eqnarray}
where

\begin{equation}\label{GRicciI}
	\left.
	\begin{aligned}
		\lambda_S&=-\frac{\mu_2F_{13}+\mu_1F_{11}}{f(2a-4br-F_1)}\\
		\lambda_g&=-\frac{4b(\mu_2F_{13}+\mu_1F_{11})}{rf(2a-4br-F_1)} 
	\end{aligned}\ \ 
	\right\rbrace	
\end{equation}
with $F_1=f_x^2+f_y^2-f(f_{xx}+f_{yy})$, $F_{13}=-2f_y^3+f_y\{4a-16br-2f_x^2+f(3f_{yy}+f_{xx}) \}-f\{-2f_xf_{xy}+f(f_{yyy}+f_{xxy}) \}$ and $F_{11}=-2f_x^3+f_x\{4a-16br-2f_y^2+f(f_{yy}+3f_{xx}) \}-f\{-2f_yf_{xy}+f(f_{xxx}+f_{xyy}) \}$.
This leads to the following:

\begin{pr}
	If $f_x^2+f_y^2-f(f_{xx}+f_{yy})-2a+4br\neq 0$, then the RT spacetime admits generalized Ricci inheritance given by
	$$\mathscr L_{\xi}S=\lambda_S S+\lambda_g g$$
	for the vector field $\xi=\mu_1 \frac{\partial}{\partial x}+\mu_2\frac{\partial}{\partial y}$, where $\lambda_S,\lambda_g$ are given in (\ref{GRicciI}).
\end{pr}

\begin{cor}
	The RT spacetime possesses Ricci inheritance given by
	$$\mathscr L_{\xi}S=\frac{\mu_2F_{13}+\mu_1F_{11}}{f(f_x^2+f_y^2-f(f_{xx}+f_{yy})-2a+4br)} S$$
	 for the vector field $\xi=\mu_1 \frac{\partial}{\partial x}+\mu_2\frac{\partial}{\partial y}$ (where $\mu_1, \mu_2$ are constants) if the relations $b=0$ and  $f_x^2+f_y^2-f(f_{xx}+f_{yy})-2a+4br\neq 0$
	are satisfied, where $F_{13}=-2f_y^3+f_y\{4a-16br-2f_x^2+f(3f_{yy}+f_{xx}) \}-f\{-2f_xf_{xy}+f(f_{yyy}+f_{xxy}) \}$ and $F_{11}=-2f_x^3+f_x\{4a-16br-2f_y^2+f(f_{yy}+3f_{xx}) \}-f\{-2f_yf_{xy}+f(f_{xxx}+f_{xyy}) \}$.
\end{cor}

\begin{cor}
	The RT spacetime possesses Ricci collineation with respect to the vector field $\xi=\mu_1 \frac{\partial}{\partial x}+\mu_2\frac{\partial}{\partial y}$ if the relations $f_x^2+f_y^2-f(f_{xx}+f_{yy})-2a+4br\neq 0$ and
		\begin{eqnarray}\label{RicciCol}
			&&\mu_2\big[-2f_y^3+f_y\{4a-16br-2f_x^2+f(3f_{yy}+f_{xx}) \}-f\{-2f_xf_{xy}+f(f_{yyy}+f_{xxy}) \}\big]\\ 
			&&+\mu_1\big[-2f_x^3+f_x\{4a-16br-2f_y^2+f(f_{yy}+3f_{xx}) \}-f\{-2f_yf_{xy}+f(f_{xxx}+f_{xyy}) \}\big]=0\notag
		\end{eqnarray}
	are satisfied, where $\mu_1, \mu_2$ are constants.
\end{cor}

The components of the Kulkarni-Nomizu products $\mathcal U=g\wedge g$, $\mathcal H= g\wedge S$ and $\mathcal D=S\wedge S$ on RT spacetime are determined as below:

$$\begin{array}{c}
	\mathcal U_{1212}= 2=\mathcal U_{2121}, 
	\mathcal U_{1221}=-2=\mathcal U_{2112},\\
	\mathcal U_{1313}=\frac{4r(d+r(-a+2br))}{f^2}=\mathcal U_{1414}=\mathcal U_{3131}=\mathcal J_{4141},\\
	\mathcal U_{1313}=-\frac{4r(d+r(-a+2br))}{f^2}=\mathcal U_{1414}=\mathcal U_{3113}=\mathcal U_{4114},\\
	
	\mathcal U_{1323}=\frac{2r^2}{f^2}=\mathcal U_{1424}=\mathcal U_{2313}=\mathcal U_{2414}=\mathcal U_{3132}	=\mathcal U_{3231}=\mathcal U_{3443}=\mathcal U_{4142}=\mathcal U_{4241}=\mathcal U_{4334},\\
	
	\mathcal U_{1323}=-\frac{2r^2}{f^2}=\mathcal U_{1442}=\mathcal U_{2331}=\mathcal U_{2441}=\mathcal U_{3123}	=\mathcal U_{3213}=\mathcal U_{3434}=\mathcal U_{4124}=\mathcal U_{4214}=\mathcal U_{4343}, 
\end{array}$$ 

$$\begin{array}{c}
	\mathcal H_{1212}=-\frac{8b}{r}=g\wedge s_{2121},\ \ 
	\mathcal H_{1221}=\frac{8b}{r}=g\wedge s_{2112},\\
	\mathcal H_{1313}=-\frac{1}{rf^2}2(d+r(-a+2br))(-2a+12br+F_1)=\mathcal H_{1414}=\mathcal H_{3131}=\mathcal H_{4141},\\
	\mathcal H_{1323}=-\frac{-2a+12br+F_1}{f^2}=\mathcal H_{1424}=\mathcal H_{2313}=\mathcal H_{2414}=\mathcal H_{3132}=\mathcal H_{3231}=\mathcal H_{4142}=\mathcal H_{4241};\\
	
	\mathcal H_{1331}=\frac{1}{rf^2}2(d+r(-a+2br))(-2a+12br+F_1)=\mathcal H_{1441}=\mathcal H_{3113}=\mathcal H_{4114};\\
	
	\mathcal H_{1332}=\frac{-2a+12br+F_1}{f^2}=\mathcal H_{1442}=\mathcal H_{2331}=\mathcal H_{2441}=\mathcal H_{3123}=\mathcal H_{3213}=\mathcal H_{4124}=\mathcal H_{4214};\\
	
	\mathcal H_{3434}=\frac{2r^2(-2a+8br+F_1)}{f^4}=\mathcal H_{4343};
	
	\mathcal H_{3443}=-\frac{1}{f^4}2r^2(-2a+8br+F_1)=\mathcal H_{4334};
\end{array}$$

$$\begin{array}{c}
	\mathcal D_{1212}=\frac{32b^2}{r^2}=\mathcal D_{2121}, \ \ 
	\mathcal D_{1221}=-\frac{32b^2}{r^2}=\mathcal D_{2112},\\
	\mathcal D_{1313}=\frac{1}{r^2f^2}16b(d+r(-a+2br))(-2a+8br+F_1)=\mathcal D_{1414}=\mathcal D_{3131}=\mathcal D_{4141}, \\   
	\mathcal D_{1323}=\frac{8b(-2a+8br+F_1)}{rf^2}=\mathcal D_{1424}=\mathcal D_{2313}=\mathcal D_{2414}=\mathcal D_{3132}=\mathcal D_{3231}=\mathcal D_{4142}=\mathcal D_{4241},\\
	\mathcal D_{1331}=-\frac{1}{r^2f^2}16b(d+r(-a+2br))(-2a+8br+F_1)=\mathcal D_{1441}=\mathcal D_{3113}=\mathcal D_{4114},\\
	\mathcal D_{1332}=-\frac{8b(-2a+8br+F_1)}{rf^2}=\mathcal D_{1442}=\mathcal D_{2331}=\mathcal D_{2441}=\mathcal D_{3123}=\mathcal D_{3213}=\mathcal D_{4124}=\mathcal D_{4214},\\
	\mathcal D_{3434}=-\frac{2(-2a+8br+F_1)}{f^4}=\mathcal D_{4343},\ \ 
	\mathcal D_{3434}=\frac{2(-2a+8br+F_1)}{f^4}=\mathcal D_{4334},
\end{array}$$ 

For the vector field $V=\frac{\partial}{\partial x}$ the components of $L_V R$ are
$$\begin{array}{c}
	(\mathscr L_V R)_{1313}=-\frac{4(d-2br^2)(d+r(-a+2br))f_x}{r^2f^3}=(\mathscr L_V R)_{1414}=(\mathscr L_V R)_{3131}=(\mathscr L_V R)_{4141},\\
	
	(\mathscr L_V R)_{1323}=-\frac{2(d-2br^2)f_x}{rf^3}=(\mathscr L_V R)_{1424}=(\mathscr L_V R)_{3132}=(\mathscr L_V R)_{3231}=(\mathscr L_V R)_{4142}=(\mathscr L_V R)_{4241},\\
	
	(\mathscr L_V R)_{1313}=\frac{4(d-2br^2)(d+r(-a+2br))f_x}{r^2f^3}=(\mathscr L_V R)_{1441}=(\mathscr L_V R)_{3113}=(\mathscr L_V R)_{4114},\\
	
	(\mathscr L_V R)_{1323}=\frac{2(d-2br^2)f_x}{rf^3}=(\mathscr L_V R)_{1442}=(\mathscr L_V R)_{3123}=(\mathscr L_V R)_{3213}=(\mathscr L_V R)_{4124}=(\mathscr L_V R)_{4214},\\

	(\mathscr L_V R)_{2313}=-\frac{2(-2b+\frac{d}{r^2})rf_x}{f^3}=(\mathscr L_V R)_{2414},\\
	
	(\mathscr L_V R)_{2331}=\frac{2(-2b+\frac{d}{r^2})rf_x}{f^3}=(\mathscr L_V R)_{2441},\\
	
	(\mathscr L_V R)_{3434}=-\frac{1}{f^5}r(4rf^3_x+f_x(8(d-ar+2br^2)+4rf^2_y-rf(3f_{yy}\\
	+5f_{xx}))+rf(-2f_yf_{xy}+f(f_{xyy}+f_{xxx})))=(\mathscr L_V R)_{4343},\\
	
	(\mathscr L_V R)_{3443}=\frac{1}{f^5}r(4rf^3_x+f_x(8(d-ar+2br^2)+4rf^2_y-rf(3f_{yy}\\
	+5f_{xx}))+rf(-2f_yf_{xy}+f(f_{xyy}+f_{xxx})))=(\mathscr L_V R)_{4334}
\end{array}$$ 
and the following relation holds:

\begin{eqnarray}
	\mathscr L_{V}R=\lambda_R R+\lambda_1 g\wedge g+\lambda_2 g\wedge S,\notag
\end{eqnarray}
where 

\begin{equation}\label{GRI}
	\left.
	\begin{aligned}
		\lambda_R&=-\frac{F_4}{f(-6d+2ar-rF_1)},\\
		\lambda_1&=-\frac{8b(-d+2br^2)f_x}{r^2f(F_1-2a+4br)}+\frac{(2ad-16bdr+8b^2r^3-dF_1)F_4}{r^3f(F_1-2a+4br)(-6d+2ar-rF_1)},\\
		\lambda_2&=-\frac{2b(-d+2br^2)f_x}{rf(F_1-2a+4br)}+\frac{(-3d+2br^2)F_4}{rf(F_1-2a+4br)(-6d+2ar-rF_1)},
	\end{aligned}\ \ 
	\right\rbrace	
\end{equation}
with $F_1=f_x^2+f_y^2-f(f_{xx}-f_{yy})$ and $F_4=(-12d+8ar-8br^2)f_x-4rf_y^2f_x+3rff_{yy}f_x-4rf_x^3+2rff_yf_{xy}-rf^2f_{xyy}+5rff_xf_{xx}-rf^2f_{xxx}$ such that $F_1-2a+4br\neq0$ and $-6d+2ar-rF_1\neq0$.
This leads to the following:

\begin{pr}
	The RT spacetime possesses generalized curvature inheritance given by 
	$$\mathscr L_{\xi}R=\lambda_R R+\lambda_1 g\wedge g+\lambda_2 g\wedge S$$
	 with respect to the vector field $\xi=\frac{\partial}{\partial x}$ if $f_x^2+f_y^2-f(f_{xx}-f_{yy})-2a+4br\neq 0$ and $r(f_x^2+f_y^2-f(f_{xx}-f_{yy}))+6d-2ar\neq 0$, where $\lambda_R, \lambda_1, \lambda_2$ are given in (\ref{GRI}),
\end{pr}

\begin{cor}
	The RT spacetime admits curvature inheritance given by
	$$\mathscr L_{\xi}R=\frac{8af_x-4f_y^2f_x+3ff_{yy}f_x-4f_x^3+2ff_yf_{xy}-f^2f_{xyy}+5ff_xf_{xx}-f^2f_{xxx}}{f(f_x^2+f_y^2-f(f_{xx}-f_{yy})-2a)} R$$
	 for the vector field $\xi=\frac{\partial}{\partial x}$ if $b=d=0$, $f_x^2+f_y^2-f(f_{xx}-f_{yy})-2a\neq 0$ and  $r(f_x^2+f_y^2-f(f_{xx}-f_{yy}))-2ar\neq 0$. 
\end{cor}

\begin{cor}
	The RT spacetime admits curvature collineation with respect to the vector field $\frac{\partial}{\partial x}$ if $b=0$, $f_x^2+f_y^2-f(f_{xx}-f_{yy})-2a\neq 0$,  $r(f_x^2+f_y^2-f(f_{xx}-f_{yy}))+6d-2ar\neq 0$ and 
	\begin{equation}\label{R-curvCol}
		(-12d+8ar)f_x-4rf_y^2f_x+3rff_{yy}f_x-4rf_x^3+2rff_yf_{xy}-rf^2f_{xyy}+5rff_xf_{xx}-rf^2f_{xxx}=0.\notag
	\end{equation}
\end{cor}

The components of the Weyl conformal (0,4) curvature tensor on RT spacetime are
$$\begin{array}{c}
	C_{1212}= -\frac{6d-2ar+rF_1}{3r^3},
	C_{1313}= \frac{1}{3r^2f^2}(d+r(-a+2br))(6d-2ar+rF_1)= C_{1414},\\
	C_{1323}=\frac{6d-2ar+rF_1}{6rf^2}=C_{1424},
	C_{3434}= \frac{r(6d-2ar+rF_1)}{3f^4}
\end{array}$$ 

Again, for the vector field $V=\frac{\partial}{\partial x}$ the components of $ \mathscr L_V C$ are
$$\begin{array}{c}
	(\mathscr L_V C)_{1212}=\frac{2(9d-2ar+rF_1)}{3r^4}=(\mathscr L_V C)_{2121},
	(\mathscr L_V C)_{1221}=\frac{-18d+4ar-2rF_1)}{3r^4}=(\mathscr L_V C)_{2112},\\
	(\mathscr L_V C)_{1313}=\frac{1}{3r^3f^2}(-4(3d^2-2adr+abr^3)+r(-d+2br^2)F_1)=(\mathscr L_V C)_{1331}\\
	=(\mathscr L_V C)_{1414}=(\mathscr L_V C)_{1441}=(\mathscr L_V C)_{3113}=(\mathscr L_V C)_{3131}=(\mathscr L_V C)_{4114}=(\mathscr L_V C)_{4141},\\
	
	(\mathscr L_V C)_{1323}=-\frac{d}{r^2f^2}=(\mathscr L_V C)_{1424}=(\mathscr L_V C)_{2313}=(\mathscr L_V C)_{2414}\\
	=(\mathscr L_V C)_{3132}=(\mathscr L_V C)_{3231}=(\mathscr L_V C)_{4142}=(\mathscr L_V C)_{4241},\\
	(\mathscr L_V C)_{1332}=\frac{d}{r^2f^2}=(\mathscr L_V C)_{1442}=(\mathscr L_V C)_{2331}=(\mathscr L_V C)_{2441}\\
	=(\mathscr L_V C)_{3123}=(\mathscr L_V C)_{3213}=(\mathscr L_V C)_{4124}=(\mathscr L_V C)_{4214},\\
	
	(\mathscr L_V C)_{3434}=\frac{6d-4ar+2rF_1}{3f^4}=(\mathscr L_V C)_{4343},
	(\mathscr L_V C)_{3443}=\frac{-6d+4ar-2rF_1}{3f^4}=(\mathscr L_V C)_{4334},
\end{array}$$

which leads to the following:
\begin{pr}
	If $f_x^2+f_y^2-f(f_{xx}-f_{yy})-2a+4br\neq 0$ and $r\{f_x^2+f_y^2-f(f_{xx}-f_{yy})\}-2ar+6d\neq 0$, then the RT spacetime possesses generalized Weyl conformal inheritance given by
	\begin{equation}
	\mathscr  L_{\xi}C=\lambda_C C+\lambda_1 g\wedge g+\lambda_2 g\wedge S\notag
	\end{equation}
	for the vector field $\xi=\frac{\partial}{\partial x}$ with
	
	\begin{equation}\label{C-curvCol}
		\left.
		\begin{aligned}
			\lambda_C&=-\frac{F_5}{f(-6d+2ar-rF_1)},\\
			\lambda_1&=-\frac{1}{{3r^3f(F_1-2a+4br)}}\big[(-12ad+4a^2r+72bdr-24abr^2)f_x\\
			&\ \ \  +(6d-4ar+12br^2)(f_y^2f_x-ff_{yy}f_x+f_x^3-ff_xf_{xx})+rf_x(f_x^4+f_y^4)\\
			&\ \ \  -2rff_y^2f_x(f_{xx}+f_{yy})+rf^2f_x(f_{xx}+f_{yy})^2+2rf_x^3(f_y^2-ff_{xx}-aff_{yy})\big],\\
			\lambda_2&=\frac{2f_x(-6d+2ar-rF_1)}{3rf(F_1-2a+4br)},	
		\end{aligned}\ \ 
		\right\rbrace	
	\end{equation}
where $F_1=f_x^2+f_y^2-f(f_{xx}-f_{yy})$ and $F_5=-12df_x+4arf_x-2rf_y^2f_x+rff_{yy}f_x-2rf_x^3+2rff_yf_{xy}-rf^2f_{xyy}+3rff_xf_{xx}-rf^2f_{xxx}$ such that $F_1-2a+4br\neq0$ and $-6d+2ar-rF_1\neq0$.
\end{pr}

\begin{cor}
	If $f_x=0$ (i.e. $f$ is a function of $y$ only) and the relations $f_y^2+ff_{yy}-2a+4br\neq 0$ and $r\{f_y^2+ff_{yy}\}-2ar+6d\neq 0$ are satisfied, then the RT spacetime possesses  Weyl conformal collineation 
	for the vector field $\frac{\partial}{\partial x}$.
\end{cor}

 The components of the concircular (0,4) curvature tensor on RT spacetime are
 $$\begin{array}{c}
 	W_{1212}=\frac{1}{6r^3}(-12d+r(-2a+12br+F_1));
 	W_{1313}=\frac{1}{3r^2f^2}(d+r(-a+2br))(6d-2ar+rF_1);\\
 	W_{1323}=\frac{6d-2ar+rF_1}{6rf^2}=\tilde{C}_{1424};
 	W_{1414}=\frac{1}{3r^2f^2}(d+r(-a+2br))(6d-2ar+rF_1);\\
 	W_{3434}=\frac{1}{6f^4}r(6d+r(-5a+6br))+5rF_1;
 \end{array}$$ 

Again, for the vector field $V=\frac{\partial}{\partial x}$ the components of $\mathscr L_V W$ are
$$\begin{array}{c}
	(\mathscr L_V W)_{1212}=\frac{1}{3r^4}(2(9d+r(a-3br))-rF_1)=(\mathscr L_V W)_{2121},\\
	(\mathscr L_V W)_{1221}=\frac{1}{3r^4}(-2(9d+r(a-3br))+rF_1)=(\mathscr L_V W)_{2112},\\
	(\mathscr L_V W)_{1313}=\frac{1}{3r^3f^2}(-4(3d^2-2adr+abr^3)+r(-d+2br^2)F_1)\\
	=(\mathscr L_V W)_{1414}=(\mathscr L_V W)_{1441}=(\mathscr L_V W)_{3131}=(\mathscr L_V W)_{4141},\\
	(\mathscr L_V W)_{1331}=\frac{1}{3r^3f^2}(4(3d^2-2adr+abr^3)-r(-d+2br^2)F_1)\\
	=(\mathscr L_V W)_{1441}=(\mathscr L_V W)_{3113}=(\mathscr L_V W)_{4114},\\
	(\mathscr L_V W)_{1323}=-\frac{d}{r^2f^2}=(\mathscr L_V W)_{1424}=(\mathscr L_V W)_{2313}\\
	=(\mathscr L_V W)_{2414}=(\mathscr L_V W)_{3132}=(\mathscr L_V W)_{3231}=(\mathscr L_V W)_{4142}=(\mathscr L_V W)_{4241},\\
	(\mathscr L_V W)_{1332}=\frac{d}{r^2f^2}=(\mathscr L_V W)_{1442}=(\mathscr L_V W)_{2331}\\
	=(\mathscr L_V W)_{2441}=(\mathscr L_V W)_{3123}=(\mathscr L_V W)_{3213}=(\mathscr L_V W)_{4124}=(\mathscr L_V W)_{4214},\\
	(\mathscr L_V W)_{3434}=\frac{1}{3f^4}(6d+2r(-5a+9br)+5rF_1)=(\mathscr L_V W)_{4343},\\
	(\mathscr L_V W)_{3443}=-\frac{1}{3f^4}(6d+2r(-5a+9br)+5rF_1)=(\mathscr L_V W)_{4334}
\end{array}$$ 
and this leads to the following:
\begin{pr}
	If $f_x^2+f_y^2-f(f_{xx}-f_{yy})-2a+4br\neq 0$ and $r\{f_x^2+f_y^2-f(f_{xx}-f_{yy})\}-2ar+6d\neq 0$, then the RT spacetime admits generalized concircular inheritance given by
	\begin{equation}
		\mathscr L_{\xi}W=\lambda_W W+\lambda_1 g\wedge g+\lambda_2 g\wedge S\notag
	\end{equation}
	 for the vector field $\xi=\frac{\partial}{\partial x}$ with
	 
	 	\begin{equation}\label{W-concirc-curvCol}
	 	\left.
	 	\begin{aligned}
	 		\lambda_W&=-\frac{F_4}{f(-6d+2ar-rF_1)},\\
	 		\lambda_1&=\frac{\big[\frac{4b}{3r^2f}(-6d+2ar-rF_1)-\frac{1}{6r^3f}(-2a+12br+F_1)\{-12d+r(-2a+12br+F_1) \}\big]F_4}{2(2a-4br-F_1)(-6d+2ar-rF_1)}\\
	 		&-\frac{(F_1-2a+12br)\{2f_yf_{xy}+f_x(f_{xx}-f_{yy})-f(f_{xyy}+f_{xxx})\}}{12r^2(2a-4br-F_1)}\\
	 		&-\frac{2b\big[2rf_x^3+f_x\{12d-4ar+2rf_y^2-rf(3f_{xx}+f_{yy})\}+rf\{-2f_yf_{xy}+f(f_{xyy}+f_{xxx}) \}\big] }{3r^2f(2a-4br-F_1)},\\
	 		\lambda_2&=\frac{f_x(-6d+2ar-rF_1)}{3rf(-2a+4br+F_1)}+\frac{(-3d+2br^2)F_4}{rf(F_1-2a+4br)(-6d+2ar-rF_1)},	
	 	\end{aligned}
	 	\right\rbrace	
	 \end{equation}
	 where $F_1=f_x^2+f_y^2-f(f_{xx}-f_{yy})$ and $F_4=(-12d+8ar-8br^2)f_x-4rf_y^2f_x+3rff_{yy}f_x-4rf_x^3+2rf_yf_{xy}-rf^2f_{xyy}+5rff_xf_{xx}-rf^2f_{xxx}$.
	 \end{pr}

\begin{cor}
	The RT spacetime possesses concircular collineation for the vector field $\frac{\partial}{\partial x}$
	if the relations $b=d=0$,   $f_x^2+f_y^2-f(f_{xx}-f_{yy})-2a\neq 0$ and $r\{f_x^2+f_y^2-f(f_{xx}-f_{yy})\}-2ar\neq 0$ are satisfied.
\end{cor}

The components of conharmonic (0,4) curvature tensor on RT spacetime are
$$\begin{array}{c}
	K_{1212}=-\frac{2(d-2br^2)}{r^3}, 
	K_{1313}=\frac{1}{r^2f^2}(d+r(-a+2br))(2(d+r(-a+4br))+rF_1)=K_{1414},\\
	K_{1323}=\frac{1}{2rf^2}(2(d+r(-a+4br))+rF_1)=K_{1424},
	K_{3434}=\frac{2r(d-2br^2)}{f^4}.
\end{array}$$ 

Again, for the vector field $V=\frac{\partial}{\partial x}$ the components of $\mathscr L_V K$ are 
$$\begin{array}{c}
	(\mathscr L_V K)_{1212}=\frac{6d-4br^2}{r^4}=(\mathscr L_V K)_{2121},
	(\mathscr L_V K)_{1221}=\frac{-6d+4br^2}{r^4}=(\mathscr L_V K)_{2112},\\
	(\mathscr L_V K)_{1313}=\frac{1}{r^3f^2}(-4(d^2-adr+br^3(3a-8br))+r(-d+2br^2)F_1)\\
	=(\mathscr L_V K)_{1414}=(\mathscr L_V K)_{3131}=(\mathscr L_V K)_{4141},\\
	(\mathscr L_V K)_{1323}=-\frac{d-4br^2}{r^2f^2}=(\mathscr L_V K)_{1424}=(\mathscr L_V K)_{2313}=(\mathscr L_V K)_{2414}=(\mathscr L_V K)_{3132}\\
	=(\mathscr L_V K)_{3231}=(\mathscr L_V K)_{4142}=(\mathscr L_V K)_{4241},\\
	(\mathscr L_V K)_{1331}=\frac{1}{r^3f^2}(4(d^2-adr+br^3(3a-8br))-r(-d+2br^2)F_1)\\
	=(\mathscr L_V K)_{1441}=(\mathscr L_V K)_{3113}=(\mathscr L_V K)_{4114},\\
	(\mathscr L_V K)_{1332}=\frac{d-4br^2}{r^2f^2}=(\mathscr L_V K)_{1442}=(\mathscr L_V K)_{2331}=(\mathscr L_V K)_{2441}=(\mathscr L_V K)_{3123}\\
	=(\mathscr L_V K)_{3213}=(\mathscr L_V K)_{4124}=(\mathscr L_V K)_{4214},\\
	
	(\mathscr L_V K)_{3434}=\frac{2(d-6br^2)}{f^4}=(\mathscr L_V K)_{4343},\ \ 
	(\mathscr L_V K)_{3443}=-\frac{2(d-6br^2)}{f^4}=(\mathscr L_V K)_{4334},
\end{array}$$ 
and the following relation holds:
\begin{eqnarray}
	\mathscr L_{V}K=\lambda_K K+\lambda_1 g\wedge g+\lambda_2 g\wedge S,\notag
\end{eqnarray}
where
\begin{equation}\label{GKI}
	\left.
	\begin{aligned}
			\lambda_K&=\frac{F_5}{f(-6d+2ar-rF_1)},\\
			\lambda_1&=\frac{16b(-d+2br^2)f_x}{r^2f(F_1-2a+4br)}-\frac{(-d+2br^2)F_5}{r^3f(-6d+2ar-rF_1)},\\
			\lambda_2&=\frac{4(-d+2br^2)f_x}{rf(F_1-2a+4br)},
	\end{aligned}\ \ 
	\right\rbrace	
\end{equation}
with $F_1=f_x^2+f_y^2-f(f_{xx}-f_{yy})$ and $F_5=(12d-4ar)f_x+2rf_y^2f_x-rff_{yy}f_x+2rf_x^3-2rff_yf_{xy}+rf^2f_{xyy}-3rff_xf_{xx}+rf^2f_{xxx}$ such that $F_1-2a+4br\neq0$ and $-6d+2ar-rF_1\neq0$.

\begin{pr}
	If $f_x^2+f_y^2-f(f_{xx}-f_{yy})-2a+4br\neq 0$ and $r\{f_x^2+f_y^2-f(f_{xx}-f_{yy})\}-2ar+6d\neq 0$, then the RT spacetime admits generalized conharmonic inheritance with respect to the vector field $\xi=\frac{\partial}{\partial x}$ given by
	$$\mathscr L_{\xi}K=\lambda_K K+\lambda_1 g\wedge g+\lambda_2 g\wedge S,$$
	 where $\lambda_K, \lambda_1,\lambda_2$ are given in (\ref{GKI}).
\end{pr}

\begin{cor}
	The RT spacetime realizes conharmonic curvature inheritance given by
	$$\mathscr L_{\xi}K= \frac{-4af_x+2f_y^2f_x-ff_{yy}f_x+2f_x^3-2ff_yf_{xy}+f^2f_{xyy}-3ff_xf_{xx}+f^2f_{xxx}}{f(2a-(f_x^2+f_y^2-f(f_{xx}-f_{yy})))} K$$
	 for the vector field $\xi=\frac{\partial}{\partial x}$ if the relations $b=d=0$, $f_x^2+f_y^2-f(f_{xx}-f_{yy})-2a\neq 0$ and $r\{f_x^2+f_y^2-f(f_{xx}-f_{yy})\}-2ar\neq 0$ hold.
\end{cor}

\begin{cor}
	The RT spacetime possesses conharmonic collineation with respect to the vector field $\frac{\partial}{\partial x}$ if $f_x=0$ (i.e., $f$ is a function of $y$ only) and the relations $f_y^2+ff_{yy}-2a+4br\neq 0$, $r\{f_y^2+ff_{yy}\}-2ar+6d\neq 0$ hold.

\end{cor}

 The components of Weyl projective (0,4) curvature tensor on RT spacetime are
$$\begin{array}{c}
	P_{1212}= -\frac{2d}{r^3}+\frac{4b}{3r},
	P_{1221}= \frac{2d}{r^3}-\frac{4b}{3r};\\
	P_{1313}=\frac{2(3d-2br^2)(d+r(-a+2br))}{3r^2f^2}, P_{1323}=\frac{3d-2br^2}{3rf^2},\\
	P_{1331}=-\frac{1}{3r^2f^2}2(d+r(-a+2br))(3d+2r(-a+br)+fF_1) = P_{1441},\\
	P_{1332}=-\frac{1}{3rf^2}(3d+2r(-a+br)+rF_1)=P_{1442}=P_{2331}=P_{2441}, \\
	P_{1414}=\frac{2(3d-2br^2)(d+r(-a+2br))}{3r^2f^2}; 
	P_{1424}= \frac{3d-2br^2}{3rf^2}=P_{2313}=P_{2414};\\
	P_{3434}=\frac{1}{3f^4}2r(3d+2r(-a+br)+rF_1)=-P_{3443}; 
\end{array}$$ 

Again, for the vector field $V=\frac{\partial}{\partial x}$, the components of $\mathscr L_V P$ are

$$\begin{array}{c}
	(\mathscr L_V P)_{1212}=\frac{6d}{r^4}-\frac{4b}{3r^2}=(\mathscr L_V P)_{2121},\ \ 
	
	(\mathscr L_V P)_{1221}=-\frac{6d}{r^4}+\frac{4b}{3r^2}=(\mathscr L_V P)_{2112},\\
	
	(\mathscr L_V P)_{1313}=\frac{-12d^2+6adr+4br^3(a-4br)}{3r^3f^2}=(\mathscr L_V P)_{1414},\\
	
	(\mathscr L_V P)_{1323}=-\frac{3d+2br^2}{3r^2f^2}=(\mathscr L_V P)_{1424}=(\mathscr L_V P)_{2313}=(\mathscr L_V P)_{2414}\\
	=(\mathscr L_V P)_{3132}=(\mathscr L_V P)_{3231}=(\mathscr L_V P)_{4142}=(\mathscr L_V P)_{4241},\\
	
	(\mathscr L_V P)_{1331}=\frac{1}{3r^3f^2}2(6d^2-5adr+2br^3(3a-4br)-r(-d+2br^2)F_1)=(\mathscr L_V P)_{1441},\\
	
	(\mathscr L_V P)_{1332}=\frac{-2b+\frac{3d}{r^2}}{3f^2}=(\mathscr L_V P)_{1442}=(\mathscr L_V P)_{2331}=(\mathscr L_V P)_{2441},\\
	
	(\mathscr L_V P)_{3123}=\frac{2b+\frac{3d}{r^2}}{3f^2}=(\mathscr L_V P)_{3213}=(\mathscr L_V P)_{4124}=(\mathscr L_V P)_{4214},\\
	
	(\mathscr L_V P)_{3131}=\frac{1}{3r^3f^2}2(-6d^2+5adr+2br^3(-3a+4br)+r(-d+2br^2)F_1)=(\mathscr L_V P)_{4141},\\
	
	(\mathscr L_V P)_{3113}=\frac{2(6d^2-3adr+2br^3(-a+4br))}{3r^3f^2}=(\mathscr L_V P)_{4114},\\
	
	(\mathscr L_V P)_{3434}=\frac{1}{3f^4}(6d+4r(-2a+3br)+4rF_1)=(\mathscr L_V P)_{4343},\\
	
	(\mathscr L_V P)_{3443}=\frac{1}{3f^4}(-6d+4r(-2a-3br)-4rF_1)=(\mathscr L_V P)_{4334}
\end{array}$$

and the following relation holds:
\begin{eqnarray}
	\mathscr L_{V}P=\lambda_P P+\lambda_1 g\wedge g+\lambda_2 g\wedge S+\lambda_3 S\wedge S,\notag
\end{eqnarray}
where

\begin{equation}\label{genGPI}
	\left.
	\begin{aligned}
			\lambda_P&=-\frac{F_9}{f(2a-4br-F_1)},\\
			\lambda_1&=\frac{F_9}{3r^5f(2a-4br-F_1)}\cdot\frac{F_7}{F_8},\\
			\lambda_2&=-\frac{F_9\{-108bdr+40b^2r^2-(9d+2br^2)(F_1-2a) \}}{3rf(F_1-2a+4br)^3} + \frac{F_{10}}{3rf(F_1-2a+4br)^2},\\
			\lambda_3&=-\frac{rF_9(-9d+2ar+2br^2-rF_1)}{3f(F_1-2a+4br)^3} - \frac{F_{12}}{3f(2a-4br-F_1)^2},
	\end{aligned}\ \ 
	\right\rbrace	
\end{equation}
with\\
\begin{equation}
		\begin{aligned}
			&F_1=f_x^2+f_y^2-f(f_{xx}+f_{yy}) \text{ such that } F_1-2a+4br\neq0,\\
			&F_7=\big[ 2(-3d+2br^2)(6br^2-2a+8br+F_1)(-2a+8br+F_1)(1+4br)\\
			&\ \ \ \ +r^2(-2a+4br+F_1) \{(-3d+2br^2)(-2a+8br+F_1)^2-16b^2(2br^2+3d-2ar+rF_1)\}  \big],\\
			&F_8=\big[(-2a+4br+F_1)\{16b^2r^2+(2a-8br-F_1)^2\}\\
			&\ \ \ \ -2(-4br+2a+F_1)(-2a+8br+F_1)(1-4br) \big]\neq0,\\
			&F_9=4af_x-8brf_x-2f^2_yf_x+ff_{yy}f_x-2f^3_x+2ff_yf_{xy}-f^2f_{xyy}+3ff_xf_{xx}-f^2f_{xxx},\\
			&F_{10}=2\big[-6adf_x+84bdrf_x-28abr^2f_x+8b^2r^3f_x+2br^2f_x\{7F_1+2f(f_{yy}-f_{xx})\}\\
			&\ \ \ \  + 4br^2f^2(f_{xxx}+f_{xyy})+3df_xF_1\big],\\
			&F_{12}=-18drf_x+8ar^2f_x-4br^3f_x-4r^2f^2_yf_x+3r^2ff_{yy}f_x-4r^2f^3_x+2r^2ff_yf_{xy}\\
			&\ \ \ \ -r^2f^2f_{xyy}+5r^2ff_xf_{xx}-r^2f^2f_{xxx}.
	\end{aligned}\notag
\end{equation}
This leads to the following:

\begin{pr}
	The RT spacetime admits generalized Weyl projective inheritance  given by
	$$	\mathscr L_{\xi}P=\lambda_P P+\lambda_1 g\wedge g+\lambda_2 g\wedge S,$$
	for the vector field $\xi=\frac{\partial}{\partial x}$,
	if the relations $f_x^2+f_y^2-f(f_{xx}+f_{yy})-2a+4br\neq 0$, $\big[(-2a+4br+F_1)\{16b^2r^2+(2a-8br-F_1)^2\}-2(-4br+2a+F_1)(-2a+8br+F_1)(1-4br) \big]\neq0$, $-9d+2ar+2br^2-r(f_x^2+f_y^2-f(f_{xx}+f_{yy}))=0$ and $-18drf_x+8ar^2f_x-4br^3f_x-4r^2f^2_yf_x+3r^2ff_{yy}f_x-4r^2f^3_x+2r^2ff_yf_{xy}-r^2f^2f_{xyy}+5r^2ff_xf_{xx}-r^2f^2f_{xxx}=0$ are satisfied with $F_1=f_x^2+f_y^2-f(f_{xx}+f_{yy})$,  where $\lambda_P,\lambda_1,\lambda_2$ are given in (\ref{genGPI}).
\end{pr}

\begin{cor}
	The RT spacetime realizes Weyl projective inheritance given by
	$$	\mathscr L_{\xi}P=\frac{f_x(4a-8br-2f^2_y+ff_{yy}-2f^2_x+3ff_{xx})+2ff_yf_{xy}-f^2f_{xyy}-f^2f_{xxx}}{f(f_x^2+f_y^2-f(f_{xx}+f_{yy})-2a+4br)} P$$
	for the vector field $\xi=\frac{\partial}{\partial x}$ if the conditions $F_1-2a+4br\neq 0$, $F_{8}\neq0$, $-108bdr+40b^2r^2-(9d+2br^2)(F_1-2a)=0$, $-9d+2ar+2br^2-rF_1$ and $F_{7}=F_{10}=F_{12}=0$ are satisfied with $F_1=f_x^2+f_y^2-f(f_{xx}+f_{yy})$, where the functions $F_{7},F_{8},F_{10},F_{12}$ are given in (\ref{genGPI}).
\end{cor}

\begin{cor}
	The RT spacetime possesses Weyl projective curvature collineation with respect to the vector field $\xi=\frac{\partial}{\partial x}$ if the relations $F_1-2a+4br\neq 0$, $F_8\neq0$ and $F_{9}=F_{10}=F_{12}=0$ hold with $F_1=f_x^2+f_y^2-f(f_{xx}+f_{yy})$, where the functions $F_{8},F_{9},F_{10},F_{12}$ are given in (\ref{genGPI}).
\end{cor}

From the above, we can state the following:

\begin{thm}
	The RT spacetime (\ref{RTM}) satisfies the following symmetry properties:
	
	\begin{enumerate}[label=(\roman*)]
		\item it admits generalized Ricci inheritance given by
		$$\mathscr L_{\xi}S=\lambda_S S+\lambda_g g$$
		for the vector field $\xi=\mu_1 \frac{\partial}{\partial x}+\mu_2\frac{\partial}{\partial y}$ ($\mu_1, \mu_2$ are constants) if $f_x^2+f_y^2-f(f_{xx}+f_{yy})-2a+4br\neq 0$, where $\lambda_S,\lambda_g$ are given in (\ref{GRicciI}),

		\item it reveals generalized curvature inheritance given by
		$$\mathscr L_{\xi}R=\lambda_R R+\lambda_1 g\wedge g+\lambda_2 g\wedge S$$
		for the vector field $\xi=\frac{\partial}{\partial x}$
		if $f_x^2+f_y^2-f(f_{xx}-f_{yy})-2a+4br\neq 0$ and $r(f_x^2+f_y^2-f(f_{xx}-f_{yy}))+6d-2ar\neq 0$, where $\lambda_R, \lambda_1, \lambda_2$ are given in (\ref{GRI}),

		\item it possesses generalized Weyl conformal inheritance given by
		$$\mathscr L_{\xi}C=\lambda_C C+\lambda_1 g\wedge g+\lambda_2 g\wedge S$$
		for the vector field $\xi=\frac{\partial}{\partial x}$ if $f_x^2+f_y^2-f(f_{xx}-f_{yy})-2a+4br\neq 0$ and $r\{f_x^2+f_y^2-f(f_{xx}-f_{yy})\}-2ar+6d\neq 0$, where $\lambda_C,\lambda_1,\lambda_2$ given in (\ref{C-curvCol}),

		\item it fulfills generalized concircular inheritance given by
		$$\mathscr L_{\xi}W=\lambda_W W+\lambda_1 g\wedge g+\lambda_2 g\wedge S$$
		for the vector field $\xi=\frac{\partial}{\partial x}$ if $f_x^2+f_y^2-f(f_{xx}-f_{yy})-2a+4br\neq 0$  and $r\{f_x^2+f_y^2-f(f_{xx}-f_{yy})\}-2ar+6d\neq 0$, where $\lambda_W, \lambda_1, \lambda_2$ are given in (\ref{W-concirc-curvCol}),

		\item it realizes generalized conharmonic inheritance given by 
		$$\mathscr L_{\xi}K=\lambda_K K+\lambda_1 g\wedge g+\lambda_2 g\wedge S$$		
		  for the vector field $\xi=\frac{\partial}{\partial x}$ if $f_x^2+f_y^2-f(f_{xx}-f_{yy})-2a+4br\neq 0$ and $r\{f_x^2+f_y^2-f(f_{xx}-f_{yy})\}-2ar+6d\neq 0$, where $\lambda_K, \lambda_1, \lambda_2$ are given in (\ref{GKI}),

		\item  it admits generalized Weyl projective inheritance given by
		 $$\mathscr L_{\xi}P=\lambda_P P+\lambda_1 g\wedge g+\lambda_2 g\wedge S$$
		 for the vector field $\xi=\frac{\partial}{\partial x}$ if the relations $f_x^2+f_y^2-f(f_{xx}+f_{yy})-2a+4br\neq 0$, $\big[(-2a+4br+F_1)\{16b^2r^2+(2a-8br-F_1)^2\}-2(-4br+2a+F_1)(-2a+8br+F_1)(1-4br) \big]\neq0$, $-9d+2ar+2br^2-r(f_x^2+f_y^2-f(f_{xx}+f_{yy}))=0$ and $-18drf_x+8ar^2f_x-4br^3f_x-4r^2f^2_yf_x+3r^2ff_{yy}f_x-4r^2f^3_x+2r^2ff_yf_{xy}-r^2f^2f_{xyy}+5r^2ff_xf_{xx}-r^2f^2f_{xxx}=0$ are satisfied with $F_1=f_x^2+f_y^2-f(f_{xx}+f_{yy})$, where $\lambda_P, \lambda_1, \lambda_2$ are given in (\ref{genGPI}).
	\end{enumerate}
\end{thm}

\begin{cor}
	The RT spacetime (\ref{RTM}) satisfies the following inheritance properties:
	
	\begin{enumerate}[label=(\roman*)]
		\item it admits Ricci inheritance given by
		$$\mathscr L_{\xi}S=\frac{\mu_2F_{13}+\mu_1F_{11}}{f(f_x^2+f_y^2-f(f_{xx}+f_{yy})-2a)} S$$
		for the vector field $\xi=\mu_1 \frac{\partial}{\partial x}+\mu_2\frac{\partial}{\partial y}$ ($\mu_1, \mu_2$ are constants) if $b=0$ and $f_x^2+f_y^2-f(f_{xx}+f_{yy})-2a\neq 0$, where the functions  $F_{13}, F_{11}$ are given in (\ref{GRicciI}),

		\item it reveals curvature inheritance given by 
		$$\mathscr L_{\xi}R=\frac{8af_x-4f_y^2f_x+3ff_{yy}f_x-4f_x^3+2ff_yf_{xy}-f^2f_{xyy}+5ff_xf_{xx}-f^2f_{xxx}}{f(f_x^2+f_y^2-f(f_{xx}-f_{yy})-2a)} R$$
		 for the vector field $\frac{\partial}{\partial x}$ if $b=d=0$, $f_x^2+f_y^2-f(f_{xx}-f_{yy})-2a\neq 0$ and $r(f_x^2+f_y^2-f(f_{xx}-f_{yy}))-2ar\neq 0$,

		\item it realizes conharmonic inheritance given by
		$$\mathscr L_{\xi}K= \frac{-4af_x+2f_y^2f_x-ff_{yy}f_x+2f_x^3-2ff_yf_{xy}+f^2f_{xyy}-3ff_xf_{xx}+f^2f_{xxx}}{f(2a-(f_x^2+f_y^2-f(f_{xx}-f_{yy})))} K$$
		for the vector field $\frac{\partial}{\partial x}$ if $b=d=0$, $f_x^2+f_y^2-f(f_{xx}-f_{yy})-2a\neq 0$ and $r\{f_x^2+f_y^2-f(f_{xx}-f_{yy})\}-2ar\neq 0$,
		
		\item  it admits Weyl projective inheritance given by
		$$	\mathscr L_{\xi}P=\frac{f_x(4a-8br-2f^2_y+ff_{yy}-2f^2_x+3ff_{xx})+2ff_yf_{xy}-f^2f_{xyy}-f^2f_{xxx}}{f(f_x^2+f_y^2-f(f_{xx}+f_{yy})-2a-+4br)} P$$
		for the vector field $\xi=\frac{\partial}{\partial x}$ if the conditions $F_1-2a+4br\neq 0$, $F_{8}\neq0$, $-108bdr+40b^2r^2-(9d+2br^2)(F_1-2a)=0$, $-9d+2ar+2br^2-rF_1$ and $F_{7}=F_{10}=F_{12}=0$ are satisfied with $F_1=f_x^2+f_y^2-f(f_{xx}+f_{yy})$, where $F_{7},F_{8},F_{10},F_{12}$ are given in (\ref{genGPI}).
	
	\end{enumerate}
\end{cor}

\begin{cor}
	The RT spacetime (\ref{RTM}) satisfies the following collineation properties:
	
	\begin{enumerate}[label=(\roman*)]
		\item it admits Ricci collineation with respect to the vector field $\mu_1 \frac{\partial}{\partial x}+\mu_2\frac{\partial}{\partial y}$ ($\mu_1, \mu_2$ are constants) if $f_x^2+f_y^2-f(f_{xx}+f_{yy})-2a+4br\neq 0$ and the relation (\ref{RicciCol}) hold,

		\item it reveals curvature collineation with respect to the vector field $\frac{\partial}{\partial x}$ if $b=0$, $f_x^2+f_y^2-f(f_{xx}-f_{yy})-2a+4br\neq 0$, $r(f_x^2+f_y^2-f(f_{xx}-f_{yy}))+6d-2ar\neq 0$ and the condition (\ref{R-curvCol})  are satisfied,

		\item it possesses Weyl conformal collineation with respect to the vector field $\frac{\partial}{\partial x}$ if  $f_x=0$ (i.e. $f$ is a function of $y$ only) and the relations $f_y^2+ff_{yy}-2a+4br\neq 0$ and $r\{f_y^2+ff_{yy}\}-2ar+6d\neq 0$ hold,

		\item it fulfills generalized concircular collineation with respect to the vector field $\frac{\partial}{\partial x}$ if the relations $b=d=0$,   $f_x^2+f_y^2-f(f_{xx}-f_{yy})-2a\neq 0$ and $r\{f_x^2+f_y^2-f(f_{xx}-f_{yy})\}-2ar\neq 0$ are satisfied,

		\item it realizes conharmonic collineation with respect to the vector field $\frac{\partial}{\partial x}$ if  $f_x=0$ (i.e., $f$ is a function of $y$ only) and the relations $f_y^2+ff_{yy}-2a+4br\neq 0$, $r\{f_y^2+ff_{yy}\}-2ar+6d\neq 0$ hold,

		\item  it admits Weyl projective collineation with respect to the vector field $\xi=\frac{\partial}{\partial x}$ if the relations $F_1-2a+4br\neq 0$, $F_8\neq0$ and $F_{9}=F_{10}=F_{12}=0$ hold with $F_1=f_x^2+f_y^2-f(f_{xx}+f_{yy})$, where the functions $F_{8},F_{9},F_{10},F_{12}$ are given in (\ref{genGPI}).
		
	\end{enumerate}
\end{cor}

\section{Conclusion}
The geometrical symmetry plays a significant role in general relativity to understand the geometry of a spacetime as it is advantageous towards the solutions of EFE and beneficial in the classification of spacetimes. The curvature collineation \cite{KLD1969} is a fundamental symmetry of spacetimes and curvature inheritance \cite{Duggal1992} is a generalized notion of curvature collineation. Indeed, during the study of curvature collineation and curvature inheritance in RT spacetime, we realize the necessity of a generalized notion of the curvature inheritance. In this paper, we have introduced the concept of \textit{generalized curvature inheritance}, which includes the notions of curvature collineation as well as curvature inheritance. Also, we have shown that the RT spacetime admits generalized curvature (resp. Ricci, Weyl conformal, concircular, conharmonic, Weyl projective) inheritance. In particular, we have proved that RT spacetime admits generalized Ricci inheritance with respect to the vector field $\mu_1 \frac{\partial}{\partial x}+\mu_2 \frac{\partial}{\partial y}$ (Proposition 4.1) and it also possesses generalized Weyl conformal (resp., concircular, conharmonic and Weyl projective) inheritance with respect to  the vector field $\frac{\partial}{\partial x}$ in Proposition 4.3 (resp., Proposition 4.4, 4.5 and 4.6). Finally, it is shown that under certain conditions the RT spacetime possesses curvature (resp. Ricci, conharmonic, Weyl projective) inheritance and (resp. Ricci, Weyl conformal, concircular, conharmonic, Weyl projective) collineation. Also, we have shown that under several conditions RT spacetime turns out to be an example of almost Ricci soliton (Proposition 3.2 and 3.3), almost $\eta$-Ricci soliton (Proposition 3.1 and 3.4), almost gradient $\eta$-Ricci soliton (Proposition 3.5) with respect to certain vector fields. \\

In this paper, certain conditions are derived in the form of partial differential equations (PDE), under which several generalized notions of symmetry and Ricci solitons are possessed by RT spacetime. But, it is extremely challenging to find their solutions as such PDE are highly non-linear (see, Theorem 3.1, 4.1 and Corollary 4.11, 4.12). Such a study is left for future investigation. Also, it is important to highlight that a number of spacetimes have not  been explored yet in the following points of views: Are they admitting any kind of generalized symmetries? What are the natures of the Ricci solitons admitted by such spacetimes? In this regard a further research is required.

\section{Appendix 
}
 
 This section is devoted to describe the programmimg code in Wolfarm Mathematica for algebraic computations of several tensors on RT spacetime, which are used to obtain the main results in this paper. First, the metric $g$ of the RT spacetime is given in (\ref{RTM}), which have been initialized for computation by the following codes:
\begin{lstlisting}[style=Matlab-editor]
    % Set dimension    
        dim = 4;
    % Input 'g' as the RT metric   
        g = {{-2(a-2*b*x2-d/x2),1,0,0} , {1,0,0,0} , {0,0,-x2^2/f^2,0} , {0,0,0,-x2^2/f^2}};
        u[1]=x1;  u[2]=x2;  u[3]=x3;  u[4]=x4;
\end{lstlisting}
Based on the metric tensor $g$, the Riemann-Christoffel symbols of kind first and second can be computed respectively as follows:
\begin{lstlisting}[style=Matlab-editor]
    % Christoffel symbols of 1st and 2nd kind:
        cs1 = Table[Simplify[D[g[[i,k]],u[j]]2+D[g[[j,k]],u[i]]2-D[g[[i,j]],u[k]]2],{i,dim},{j,dim},{k,dim}];
        cs2 = Table[Simplify[Sum[G[[h,k]] cs1[[i,j,k]],{k,1,dim}]],{i,dim},{j,dim},{h,dim}];
\end{lstlisting}
Using the Riemann-Christoffel symbols, the Riemann tensor of type (1,3) and 
the Riemann tensor of type (0,4) are computed respectively by the following code: 
\begin{lstlisting}[style=Matlab-editor]
    % Compute 'rr' as the (1,3)-type curvature tensor:
        For[i=1,i<=dim,i++,
         For[j=1,j<=dim,j++,
          For[k=1,k<=dim,k++,
           For[h=1,h<=dim,h++,
            rr[h,i,j,k] = FullSimplify[D[cs2[[k,i,h]],u[j]]-D[cs2[[j,i,h]],u[k]] + Sum[cs2[[k,i,l]] cs2[[j,l,h]],{l,1,dim}] - Sum[cs2[[j,i,l]] cs2[[k,l,h]],{l,1,dim}]]]]]];
            
    % Compute 'RR' as the (0,4)-type curvature tensor:
        RR=Table[Sum[g[[h,q]] rr[h,i,j,k],{h,1,dim}],{q,dim},{i,dim},{j,dim},{k,dim}];
\end{lstlisting}
The Ricci tensor $S$ and scalar curvature $\kappa$ can be computed respectively by the following code:
\begin{lstlisting}[style=Matlab-editor]
    % Calculate 'ss' as the Ricci tensor:
            ss = Table[FullSimplify[Sum[G[[l,k]] RR[[l,i,j,k]], {l,1,dim}, {k,1,dim}]], {i,dim}, {j,dim}];
    % Compute 'kappa' as the scalar curvature:
            kappa = FullSimplify[Sum[G[[i,j]] ss[[i,j]], {i,dim}, {j,dim}]];
\end{lstlisting}
By using the above codes we have determined the (0,4)-type Riemann curvature tensor $R$, Ricci tensor $S$ and the scalar curvature $\kappa$ of the RT spacetime, which is shown in Section 2. Next, the Kulkarni-Nobizu product tensors $g\wedge g$, $g\wedge S$, $S\wedge S$ can be computed respectively by the following code:
\begin{lstlisting}[style=Matlab-editor]
        gag = Table[2 FullSimplify[g[[i,l]]*g[[j,k]] - g[[i,k]]*g[[j,l]]], {i,dim}, {j,dim}, {k,dim}, {l,dim}];
\end{lstlisting}

\begin{lstlisting}[style=Matlab-editor]
        gas = Table[ FullSimplify[ g[[i,l]]*ss[[j,k]] + g[[j,k]]*ss[[i,l]] - g[[i,k]]*ss[[j,l]] - g[[j,l]]*ss[[i,k]]], {i,dim}, {j,dim}, {k,dim}, {l,dim}]
\end{lstlisting}

\begin{lstlisting}[style=Matlab-editor]
        sas = Table[ 2 FullSimplify[ss[[i,l]]*ss[[j,k]] - ss[[i,k]]*ss[[j,l]]], {i,dim}, {j,dim}, {k,dim}, {l,dim}];
\end{lstlisting}
By employing these codes, the components of $g\wedge g$, $g\wedge S$ and $S\wedge S$ of RT spacetime are obtained as it is exhibited in Section 4. Next, for facilitating further computations of the several necessary tensors we define `Fa', `tenp1', `Li' in the following codes:
\begin{lstlisting}[style=Matlab-editor]
        Fa = A @@@ Transpose[{Table[i, {i,dim}]}];
        tenp1[A_, R_] := Table[R[i]*A[m], {i,dim}, {m,dim}]
        Li[v1_, v2_] := Table[FullSimplify[Sum[v1[[j]] D[v2[[i]], u[j]], {j,1,dim}] - Sum[v2[[j]] D[v1[[i]], u[j]], {j,1,dim}]], {i,dim}]
\end{lstlisting}
Now, if $V$ is vector field and $X$ is a tensor of type (0,2), the Lie derivative ($\mathscr L_V X$) of $X$ with respect to $V$ can be computed as follows:
\begin{lstlisting}[style=Matlab-editor]
    % 'Lie2' is the Lie derivative of (0,2)-type tensors
        Lie2[v_][X_] := Table[FullSimplify[Sum[v[[k]] D[X[[i,j]], u[k]], {k,1,dim}] - Sum[Li[v,e[i]][[k]] X[[k,j]], {k,1,dim}] - Sum[Li[v,e[j]][[k]] X[[i,k]], {k,1,dim}]], {i,dim}, {j,dim}]
\end{lstlisting}
Again, for a (0,4)-type tensor $Y$, the Lie derivative ($\mathscr L_V Y$) of $Y$ with respect to vector field $V$ can be calculated by the following codes:
\begin{lstlisting}[style=Matlab-editor]
    % 'Lie4' is the Lie derivative of (0,2)-type tensors
        Lie4[v_][Y_] := Table[FullSimplify[ Sum[v[[k]] D[Y[[i,j,m,n]], u[k]], {k,1,dim}] - Sum[Li[v,e[i]][[k]] Y[[k,j,m,n]], {k,1,dim}] - Sum[Li[v,e[j]][[k]] Y[[i,k,m,n]], {k,1,dim}] - Sum[Li[v,e[m]][[k]] Y[[i,j,k,n]], {k,1,dim}] - Sum[Li[v,e[n]][[k]] Y[[i,j,m,k]], {k,1,dim}]], {i,dim}, {j,dim}, {m,dim}, {n,dim}]
\end{lstlisting}
From the above tensor calculations, the existence of generalized curvature inheritance (Proposition 4.2) on the RT spacetime is shown by utilizing the following codes:
\begin{lstlisting}[style=Matlab-editor]
        Solve[Lie4[v][RR] == x*RR + y*gag + z*gas + w*sas, {x,y,z,w}]
\end{lstlisting}
Similarly, we can prove that the RT spacetime admits generalized Ricci (resp. Weyl conformal, concircular, conharmonic, Weyl projective) inheritance, which is shown in Proposition 4.1 - 4.6. Also, the nature of Ricci soliton (Proposition 3.1 - 3.4 and Theorem 3.1) admitted by RT spacetime can be figured out by the following codes: 
\begin{lstlisting}[style=Matlab-editor]
        Solve[Lie2[v][g] + x*ss == y*g + z*tenp1[A,A], {x,y,z} \[Union] Fa]
\end{lstlisting}
 The readers are requested to go through \cite{Shaikh_Kundu_2015} for the detailed classification of semi-Riemannian spaces by using the algebraic computations in Wolfram Mathematica regarding various tensors like conformal, conharmonic, concircular, Weyl projective curvature tensor.\\

\section*{Acknowledgment}
The second author is grateful to The Council of Scientific and Industrial Research (CSIR File No.: 09/025(0253)/2018-EMR-I), Govt. of India, for the award of SRF (Senior Research Fellowship). All the algebraic computations of Section 3 and 4 are performed by a program in Wolfram Mathematica developed by the first author A.A.Shaikh.\\

\end{document}